\documentclass[letterpaper, 10 pt, journal, twoside]{ieeetran}

\IEEEoverridecommandlockouts                              

\usepackage{tikz}
\usepackage{graphics} 
\usepackage{epsfig} 
\usepackage{amsmath} 
\usepackage{amssymb}  
\usepackage{url,algorithm}
\usepackage{multirow}
\usepackage{color}

\usepackage{enumitem}
\renewcommand{\theenumi}{\arabic{enumi}}

\renewcommand{\theenumii}{\arabic{enumii}}

\renewcommand{\theenumiii}{\arabic{enumiii}}

\renewcommand{\theenumiv}{\arabic{enumiv}}

\newcommand{\smallmat}[1]{\left[ \begin{smallmatrix}#1 \end{smallmatrix} \right]}

{\normalfont}{\normalfont}

\newcommand{\TT}{T}
\newcommand{\Ts}{T_s}
\newcommand{\TMPC}{T_{\text{MPC}}}
\newcommand{\nin}{n_u}
\newcommand{\nout}{n_y}
\newcommand{\rK}{g}
\newcommand{\JC}{\tilde{J}}

\newcommand{\Nu}{N_u}
\newcommand{\Nc}{\Nu}
\newcommand{\Np}{N_p}
\newcommand{\Kin}{\mathcal{K}}
\newcommand{\Min}{\mathcal{M}}
\newcommand{\Nin}{N_{\text{in}}}

\newcommand{\matrice}[2]{\left[\hspace*{-.1cm}\begin{array}{#1} #2 \end{array}\hspace*{-.1cm}\right]}

\title{
Performance-oriented model learning for data-driven MPC design
}

\author{Dario Piga, Marco Forgione, Simone Formentin, Alberto Bemporad
	\thanks{This work was partially supported by the H2020-723248 project
		\emph{DAEDALUS - Distributed control and simulation platform to support an
		ecosystem of digital automation developers} and by the Lombardia region and the
		Cariplo foundation, under the project \emph{Learning to Control (L2C)}, no.
		2017-1520.}
		\thanks{D. Piga and M. Forgione are with IDSIA Dalle Molle Institute for Artificial Intelligence SUPSI-USI, Manno, Switzerland. {\tt\small dario.piga@supsi.ch; marco.forgione@supsi.ch} }
	\thanks{S. Formentin is with the Dipartimento di Elettronica e Informazione, Politecnico di Milano, Milano, Italy. {\tt\small simone.formentin@polimi.it}}  
	\thanks{A. Bemporad is with IMT School for Advanced Studies Lucca,  Lucca, Italy 
		{\tt\small alberto.bemporad@imtlucca.it}}   %
}

\begin{document}

\maketitle
\thispagestyle{empty}

\begin{abstract}
Model Predictive Control (MPC) is an enabling technology in applications requiring controlling physical processes in an optimized way under constraints on inputs and outputs. However, in MPC closed-loop performance is pushed to the limits only if the plant under control is accurately modeled; otherwise, robust architectures need to be employed, at the price of reduced performance
due to worst-case conservative assumptions. In this paper, instead of adapting the controller to handle uncertainty, we adapt the learning procedure so that the prediction model is selected to provide the best closed-loop  performance. More specifically, we apply for the first time the above ``identification for control'' rationale to hierarchical MPC using data-driven methods and Bayesian optimization. 
\end{abstract}


\section{Introduction}
\IEEEPARstart{N}{owadays}, \textit{Model Predictive Control} (MPC) has become the most popular advanced control technology for several complex engineering applications \cite{borrelli2017predictive}. Apart from computational aspects, it is widely recognized that one key practical challenge in MPC arises when dealing with uncertainty, especially when the prediction model is identified using open-loop data taken from a specific operation of the plant~\cite{mesbah2016stochastic}.

In case of partially known systems, traditional MPC approaches exhibit some degree of robustness, so that marginal robust performance can be guaranteed. When intrinsic robustness of deterministic
MPC is not enough, robust MPC \cite{falugi2014getting} and stochastic MPC  \cite{mesbah2016stochastic} approaches have been developed to take into account  uncertainties. However, regardless of the specific technique, increasing robustness of the MPC controller usually leads to   conservative performance \cite{bemporad1999robust}.

While there is usually a separation between model identification and control design, an alternative approach to managing uncertainty in designing control systems is to revisit the identification process as a procedure \textit{to be designed} by bearing the final control application in mind. Such a rationale is known as \textit{Identification for Control} (I4C) and has been widely studied for fixed-order (oftentimes, PID) control of \emph{Linear Time-Invariant} (LTI) systems \cite{gevers2005identification}. According to I4C, the best model for control may \textit{not} be the one providing the least output prediction errors, but the one providing the best performance on the true system when in closed loop with the associated model-based controller.

As far as we are aware of, the I4C modeling approach has \textit{never} been applied to MPC control. Learning techniques have instead been shown to be useful for iterative MPC tasks in \cite{rosolia2018learning} and in reinforcement learning applied to MPC~\cite{ZGB19}. Furthermore, data-driven approaches have been proposed for direct MPC optimization using open- and closed-loop data, see, e.g., \cite{kadali2003data,PiFoBe2017}. Although the above approaches are powerful tools for control design in case of unknown systems, they fail to provide a mathematical (albeit control-oriented) description of the plant. Indeed, the latter can often be   useful for physical interpretation, performance monitoring, and diagnosis \cite{schafer2004multivariable}.

In this work, we propose an \textit{Identification for (Model-Predictive) Control} approach aimed at finding the best prediction model for MPC  from experimental data, by considering the control objective directly in the model learning phase. We propose a  hierarchical architecture, typically employed in several industrial applications, in which the inner controller is a parametric filter (e.g., a PID controller) aimed to stabilize the system at a fast pace, whereas the outer loop plays the role of a reference governor (RG)~\cite{garone2017reference,BCM97} with a twofold goal: ($i$) boosting the performance of the inner loop and ($ii$) handling the signals constraints due, e.g., to actuator bounds or system limitations. Within this framework, the RG is typically an MPC law based on a model of the inner loop. According to the I4C philosophy, we propose a change of perspective and treat such a model as a \textit{design parameter} instead. Such a parameter will be iteratively optimized, together with the inner controller, using closed-loop data collected on the plant and Bayesian optimization. Finally, we show that, using the same rationale and tools, also 
the prediction horizon, a critical parameter to tune in MPC, can be optimized from data.

For the sake of completeness, the first use of Bayesian optimization in control-oriented identification was proposed in \cite{bansal2017goal}, based on a simpler control scheme. The same hierarchical architecture was instead addressed in \cite{PiFoBe2017} to design the controller from data, but without providing an MPC-oriented model of the plant.

The remainder of the paper is as follows. In Section \ref{sec:problem} the control problem of interest is formally stated. The hierarchical architecture is introduced in Section \ref{sec:architecture}, where also the parameterization of each block is described (and motivated) in detail. The proposed strategy is described in Section \ref{sec:main}, where a discussion on how to practically restrict the parameter space is also provided. Section \ref{sec:example} illustrates the performance of the method on a benchmark   example. 


\section{Problem formulation}\label{sec:problem}

Consider a \emph{multi-input multi-output} (MIMO) plant $\mathcal{S}$, with input $u\in \mathbb{R}^{\nin}$ and output $y\in \mathbb{R}^{\nout}$ signals sampled at a regular time interval $\Ts$. We aim at synthesizing a controller $\mathcal{C}$ for $\mathcal{S}$  
such that the controlled closed-loop system achieves a desired engineering objective defined in terms of minimization of a cost $J(y_{1:\TT},u_{1:\TT})$, where $y_{1:\TT}$ (resp. $u_{1:\TT}$) denotes the sequence of output  (resp. input) signals measured at time steps $t=1,\ldots,\TT$, and $\TT$ is the length (measured in number of samples) of the experiment where the closed-loop performance is measured. Besides minimizing the cost $J(y_{1:\TT},u_{1:\TT})$,  
the following constraints on inputs and outputs should be satisfied:
\begin{subequations} \label{eq:constraints}
	\begin{align}
	u_{\mathrm{min}} &\leq u(t) \leq u_{\mathrm{max}}, \label{eq:constraintsJa}\\ 
	\Delta u_{\mathrm{min}} &\leq u(t)-u(t-1)  \leq \Delta u_{\mathrm{max}}, \label{eq:constraintsJb}\\ 
	y_{\mathrm{min}} &\leq y(t)  \leq y_{\mathrm{max}}, \quad  t=1,\ldots,\TT. \label{eq:constraintsJc} 
	\end{align}
\end{subequations}
Constraints~\eqref{eq:constraints} are generally imposed by actuator limitations or might reflect safety conditions. The control design problem is formulated as the following optimization problem:
\begin{align} \label{eqn:costCLv1}
\min_{\mathcal{C} \in \mathcal{\mathbf{C}}}& \quad J(y_{1:\TT},u_{1:\TT}) \quad 
\textrm{s.t.}   \quad \text{(\ref{eq:constraints})},
\end{align}
with $\mathcal{\mathbf{C}}$ denoting the set of controller candidates.

We rewrite constraints~(\ref{eq:constraints}) as $h_i(t) \geq 0$, $i=1,\ldots,6$,
with
\begin{subequations} \label{eqn:poscons}
	\begin{align}
	&  h_1(t)= u(t) - u_{\mathrm{min}}, \ \ \ \ \  h_2(t)=   u_{\mathrm{max}}- u(t),\\
	&  h_3(t)= u(t)-u(t-1) - \Delta u_{\mathrm{min}} \geq 0,\\
	&  h_4(t)= \Delta u_{\mathrm{max}} -  u(t)+u(t-1)\geq 0 ,\\
	&  h_5(t)= y(t) - y_{\mathrm{min}}, \ \ \  h_6(t)=   y_{\mathrm{max}}- y(t),
	\end{align}
\end{subequations} 
and treat them with penalty functions 
\begin{align} \label{eqn:costCL}
\min_{\mathcal{C} \in \mathcal{\mathbf{C}}}& \quad \JC(y_{1:\TT},u_{1:\TT},\mathcal{C} )
\end{align}
where
\begin{align}
\JC(y_{1:\TT},u_{1:\TT}) = J(y_{1:\TT},u_{1:\TT}) + \sum_{t=1}^{\TT}\sum_{i=1}^6b_t(h_i(t))
\label{eq:finalJC}
\end{align}
and $b_t:\mathbb{R}\to \mathbb{R}$ are (possibly time-varying) barrier functions. 

Assuming zero initial conditions, clearly $y_{1:\TT}$, $u_{1:\TT}$ in~\eqref{eq:finalJC} are only functions of the controller $\mathcal{C}$ and of the process model $\mathcal{S}$. Rather than first fixing a model for $\mathcal{S}$ (either from  first-principle physical laws or using system identification techniques), we follow a \emph{performance-driven} control design paradigm
and leave the model of $\mathcal{S}$ as a degree of freedom, used to
minimize the closed-loop cost $\JC(y_{1:\TT},u_{1:\TT})$. 


\section{Control architecture}\label{sec:architecture}

\begin{figure}[!b]
	\begin{tikzpicture}[scale=.77]
	\node [draw, rectangle,minimum width=1.2cm,minimum height=1cm,thick] (Gp)  at (4,1.5) {$\mathcal{S}$};
	\node [above] (yo)  at (5.2,1.5) {};
	\node [draw, rectangle,minimum width=1.2cm,minimum height=1cm, thick] (Kp)  at (1.5,1.5) {$\Kin$};
	\draw[->] [thick]   (Kp) -- (Gp);
	\node [] (plus)  at (6.9,1.5) {};
	\node [draw, circle,scale=0.5] (pluserr)  at (0,1.5) {};
	\node [above right] (r)  at (-1.5,1.5) {\footnotesize $g(t)$};
	\node [below left,scale=0.5]    at (0,1.35) {$-$};
	\draw[->] [thick] (Gp) -- (plus);
	\node[above] (ynode22) at  (6.5,1.5) {\footnotesize  $y(t)$};
	\draw[->] [thick]  (6.0,1.5) -- (6.0,-0.8) -- (-3.3,-0.8) -- (-3.3, 1.2) -- (-2.95,1.2);
	\draw[->] [thick] (pluserr) -- (Kp);
	\node [above] (u)  at (2.75,1.5) {\footnotesize $u(t)$};
	\node [above] (e)  at (0.4,1.5) {\footnotesize $e(t)$};
	\draw[<-] [thick] (pluserr) -- (0,0.5) -- (5.3,0.5) -- (5.3,1.5);
	\draw[dashed] (-0.65,-0.3) rectangle(5.6,2.5);
	\node [above right] (M)  at (-0.7,-0.25) {$\Min$};
	\node [draw, rectangle,minimum width=1.2cm,minimum height=1cm, thick] (MPC)  at (-2.2,1.5) {$\mathrm{MPC}$};
	\draw[->] [thick] (MPC) -- (pluserr);
	\node [left] (rout)  at (-3.8,1.5) {};
	\node [above] (rout2)  at (-3.5,1.5) {\footnotesize $  r(t)$};
	\draw[->] [thick] (rout)--(MPC); 
	\end{tikzpicture}   
	\caption{Proposed hierarchical control architecture. $\mathcal{S}$: plant to be controlled; $\Kin$: inner  controller; $\Min$: inner closed-loop system; $r(t)$:   reference to be tracked.}  \label{fig:hierarchical_outinner}  \vspace*{-0.0cm}
\end{figure}
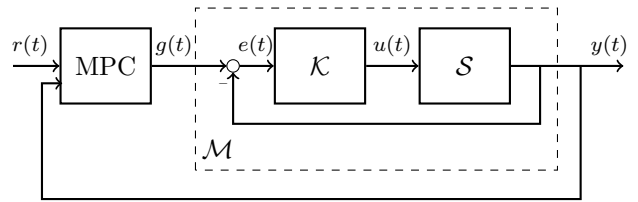

We adopt the hierarchical, multi-rate, reference-governor control architecture  in  Fig.~\ref{fig:hierarchical_outinner}, consisting of: 
\begin{itemize} 
	\item an inner low-level  controller   $\Kin$
	which operates at sampling time $\Ts$ and  it is
	mainly used to handle fast
	dynamics of the system. This controller introduces a
	degree-of-freedom in the control design and, in case
	of unstable plants $\mathcal{S}$, it might also stabilize the inner
	closed-loop system $\mathcal{M}$. Nevertheless, the latter is not a
	required condition in our design approach.
 \item an outer MPC to enhance performance of the inner loop $\Min$ an to enforce constraints~\eqref{eq:constraintsJc}. The MPC operates at a sampling time $\TMPC$ that is an integer multiple of $\Ts$, \emph{i.e.}, $\TMPC=N\Ts$ with $N   \in\mathbb{N}$. 
Setting $\TMPC$ larger than $\Ts$ (thus, $N > 1$) may be needed  to solve the constrained optimization problem on line, \emph{i.e.}, within the MPC sampling time $\TMPC$.
\end{itemize}
In standard RG approaches, the outer MPC requires a prediction model of the inner loop $\Min$. 
In accordance with the \emph{performance-driven}  approach proposed in this paper, 
we treat such a model as a \textit{design parameter}
and look for the model providing the best closed-loop performance according to the performance index $\JC(y_{1:\TT},u_{1:\TT})$. 
In particular, as detailed in the following, a model of the plant $\mathcal{S}$ will be used neither to design the controllers  nor to evaluate  the  performance index $\JC(y_{1:\TT},u_{1:\TT})$, which will be instead measured directly from closed-loop experiments performed on the actual plant.

\subsection{Inner controller parameterization} \label{Sec:inner}
The inner controller $\Kin$ is parameterized by a vector $\theta\in\mathbb{R}^{n_\theta}$. For instance, $\Kin$ can be a simple discrete-time \emph{proportional-integral-derivative} (PID) controller, with sampling time $\Ts$ and discrete-time transfer function
\begin{equation}
\label{eq:PID}
\Kin(z,\theta) = \theta_P + \theta_I\Ts\frac{1}{z-1} + \theta_D \frac{N_d}{1+N_d\Ts\frac{1}{z-1}}, 
\end{equation}
where $\theta = \left[ \theta_P \ \  \theta_I\ \ \theta_D \right]'$ is the design parameter vector   and  $N_d \gg 1$  limits the high-frequency gain of the PID controller. Although $N_d$ may  be treated as a design parameter, its tuning  is generally not critical and thus   not included in $\theta$.   

\subsection{Outer MPC parameterization} \label{Sec:outerMPC}
The most important component of the outer MPC is the model 
used to predict the output $y$ and input $u$ as a function of the MPC command $\rK$. %
Let $M$ be the dynamical model from $\rK$ to $\smallmat{y\\u}$, described in the  state-space representation \begin{equation}  
\left\{\begin{array}{rcl}
\xi(t+1)&=&A_{M}\xi(t)+B_{M}\rK(t)\\
\matrice{c}{y(t)\\u(t)}&=&C_{M}\xi(t)+D_{M}\rK(t),
\end{array}
\right.
\label{eqn:SITO}
\end{equation}
where $\xi\in\mathcal{R}^{n_\xi}$ is the state of the closed-loop model.
For instance, in the case of a single-input-single-output  plant, 
the $2\times1$ transfer matrix $M$ can be modelled as a pair of 
transfer functions with the same  poles.
Let $\mu\in\mathbb{R}^{n_\mu}$ be the vector obtained by stacking the entries of  $A_{M}, B_{M}, C_{M}$ and $D_{M}$.

At each time instant $t$ integer multiple of the MPC sampling time $\TMPC$ (\emph{i.e.}, $t=h\TMPC$ with $h \in \mathbb{N}$), the outer MPC solves the  minimization problem
\begin{subequations} \label{eq:MPC}
	\begin{align}
	&  \nonumber \min_{\left\{\rK(t+k|t)\right\}_{k=1}^{\Nu}, \epsilon} Q_y \!\sum_{k=1}^{\Np} \left(y(t+k|t)-r(t+k)\right)^2+\\ 
	&  \nonumber \quad\quad\quad\quad+ Q_u\! \sum_{k=1}^{\Np} \left(u(t+k|t)-u_{\rm ref}(t+k)\right)^2 \!+\!\\
	&   \quad\quad\quad\quad + Q_{\Delta u}\!\sum_{k=1}^{\Np}\!\!\left(\!u(t+k|k)\!-\!u(t+k-1|t)\right)^2   \!\!+\! Q_\epsilon \epsilon^2\\
	&  \mathrm{s.t. \ }    \matrice{c}{y(t+k|t)\\u(t+k|t)} = M\left(\mu,
	\rK(t+k|t)\right),  \ \   k=1,\ldots,\Np  \\
	& \quad y_{\mathrm{min}} \!-\! V_y\epsilon\leq y(t+k|t) \leq  y_{\mathrm{max}}+V_y\epsilon, \  k=1,\ldots,\Np  \\
	& \quad   u_{\mathrm{min}}\!-\!V_u\epsilon \leq u(t+k|t) \leq  u_{\mathrm{max}}+V_u\epsilon, \  k=1,\ldots,\Np \\
	& \quad  \Delta u_{\mathrm{min}} -V_{\Delta u}\epsilon \leq \Delta u(t+k|t), \ \    k=1,\ldots,\Np \label{eq:MPC:Ducons1}  \\
	& \quad    \Delta u(t+k|t) \leq  \Delta u_{\mathrm{max}}+V_{\Delta u}\epsilon, \ \  k=1,\ldots,\Np \label{eq:MPC:Ducons2}\\ 
	& \quad  \rK(t+\Nu+j|t)= \rK(t+\Nu|t),  \ \  j=1,\ldots, \Np-\Nu	 \label{eqn:MPCconst}
	\end{align}
\end{subequations}
where $\Delta u(t+k|t) =  u(t+k|t)- u(t+k-1|t)$,   $\Np$ and $\Nu$ are the prediction and control horizon, respectively, $Q_{y}$, $Q_{u}$, $Q_{\Delta u}$, $Q_\epsilon$ are nonnegative weights, $u_{\rm ref}$ and $r$ are the input and output references, respectively,  $V_y$, $V_u$, $V_{\Delta u}$ are positive vectors that are used to soften the  constraints on plant's input and output, so that~\eqref{eq:MPC} always admits a solution. According to standard MPC design, in case $\Nu<\Np$,  the constraint \eqref{eqn:MPCconst} enforces a constant value of   $g$ from time $\Nu$ to $\Np$.  
The reader is referred to \cite{borrelli2017predictive} for an overview on  MPC design.

We can also treat the
prediction horizon $\Np$ as a design parameter, and denote by 
$\nu=\left[\mu'\ \ \Np \right]'$, $\nu\in\mathbb{R}^{n_\mu}\times\mathbb{N}$ the overall vector of tuning parameters. The control horizon $\Nu$ determines, together with the number of constraints in~\eqref{eq:MPC}, the computational complexity of the outer MPC controller. Therefore, it is usually fixed by the available online throughput. Alternatively, we can set $\Nu=\Np$.

The remaining MPC parameters ($\Nc$, $Q_{y}$, $Q_{u}$, $Q_{\Delta u}$, $Q_\epsilon$, $V_y$, $V_u$ and $V_{\Delta u}$) are   treated as a specification
of the desired closed-loop performance, and therefore not optimized.
More generally, we could decouple the MPC quadratic cost in~\eqref{eq:MPC} from
the closed-loop performance index $\JC(y_{1:\TT},u_{1:\TT})$. For instance, $\JC(y_{1:\TT},u_{1:\TT})$ can be a general, possibly non-convex function 
reflecting engineering or economic goals,  while the cost of the MPC~\eqref{eq:MPC} is  quadratic to facilitate online optimization. Indeed, in case the augmented model $M$ is LTI as in~\eqref{eqn:SITO},  problem~\eqref{eq:MPC} reduces to a  \emph{quadratic programming} (QP) problem whose solution can be   computed both \emph{offline} using multiparametric quadratic programming~\cite{BEMPORAD20023} or online using dedicated QP solvers based, \emph{e.g.}, on interior-point algorithms~\cite{WaBo10}, fast gradient projection~\cite{PaBe14}, or active set methods~\cite{ferreau2008online}. 


\section{Performance-driven parameter tuning}\label{sec:main}

Based on the controller parametrization introduced in the previous section, the closed-loop performance cost $\JC(y_{1:\TT},u_{1:\TT})$ is  a function of vectors $\theta$ and $\nu$ parametrizing the inner controller $\Kin$ and the outer MPC, respectively. Thus, under the hierarchical architecture of Fig.~\ref{fig:hierarchical_outinner}, the original control design problem~\eqref{eqn:costCL} is equivalent to
\begin{align} \label{eqn:costCLtot}
\min_{\theta,\nu}& \quad \JC(y_{1:\TT},u_{1:\TT};\theta,\nu).
\end{align}

\subsection{Bayesian optimization for parameter selection}

\begin{algorithm}[!b]
	\caption{Bayesian optimization for controller design}
	\label{algo:BO}
	\begin{enumerate}[label=\arabic*., ref=\theenumi{}]
		\item \label{algo:BOv1} \textbf{initialize} parameters \emph{vs} performance set
		\begin{align*} \mathcal{D} \leftarrow \{(\theta_{1:\Nin},\nu_{1:\Nin}),\JC_{1:\Nin}\}; 
		\end{align*} 
		\item  \textbf{for} $i=\Nin,\ldots,i_{\mathrm{max}}-1$ \textbf{do}
		\begin{enumerate}[label=\theenumi{}.\arabic*., ref=\theenumi{}.\theenumii{}]
			\item \label{algo:BOv2} based on the data  $\mathcal{D}$,  \textbf{train} a GP  approximating $\JC$;
			\item \label{algo:BOv2.2} based on the GP,  \textbf{define} acquisition function $\alpha(\theta,\nu|\mathcal{D})$;
			\item \textbf{compute} next controller parameters $\theta_{i},\nu_{i}$ as  
			\begin{align*}
			\theta_{i+1},\nu_{i+1} \leftarrow \arg\max\limits_{\theta,\nu}\alpha(\theta,\nu|\mathcal{D});
			\end{align*}
			\item \label{algo:BOv3} \textbf{perform} closed-loop experiment and \textbf{measure} performance index $\JC_{i+1}$;
			\item \textbf{augment} the training set {\small $ \mathcal{D}\! \leftarrow \! \mathcal{D} \cup \{(\theta_{i+1},\nu_{i+1}),\JC_{i+1}\}$};
			\item   \textbf{exit for loop} if a termination criterion is met;
		\end{enumerate}
		\item \textbf{end for}
		\item \textbf{compute} optimal  parameters as $\theta_{i^{\star}}$ and $\nu_{i^{\star}}$, with 
		\begin{align*}
		i^{\star}=\arg \max\limits_{i}\JC_{i};
		\end{align*}
	\end{enumerate}
	\vspace*{-.0cm}\hrule\vspace*{.0cm}
	~~\textbf{Output}:  Optimal controller parameters $\theta_{i^{\star}}$ and $\nu_{i^{\star}}$. 
\end{algorithm}

The design problem~\eqref{eqn:costCLtot} is solved through the \emph{Bayesian optimization} (BO) strategy \cite{shahriari2016taking} outlined in Algorithm~\ref{algo:BO}. The algorithm is initialized (step~\ref{algo:BOv1}) by  performing $\Nin \geq 1$ closed-loop experiments for $\Nin$ different (e.g., randomly chosen) values of controller parameters $\theta_i$ and $\nu_i$ (with $i=1,\ldots,\Nin$). For each pair $(\theta_i,\nu_i)$, a closed-loop experiment is performed and the performance index $\JC_i$ is measured. In this way,  an initial set $\mathcal{D} = \{(\theta_{1:\Nin},\nu_{1:\Nin}),\JC_{1:\Nin}\}$ of parameters and corresponding performance $\JC$ is constructed, with $\theta_{1:\Nin}$ denoting the sequence $\theta_{i}$ for $i=1,\ldots,\Nin$. In practice, the experiment can be interrupted and large cost assigned to $\JC_i$  in case of safety constraint violations.


The algorithm is then iterated until a stopping criterion is met (e.g., maximum number of iterations reached). At each iteration $i \geq \Nin$, the following two steps are performed:
\begin{itemize}
	\item \emph{Learning phase} (Step~\ref{algo:BOv2}). In this step, a Gaussian Process (GP) describing our ``best guess'' of the cost $\JC(y_{1:\TT},u_{1:\TT};\theta,\nu)$ 
corresponding to the design parameters $\theta$ and $\nu$ is fitted on the available data  $\mathcal{D}$. 
	
	Under the prior assumption that the   cost $\JC$ is generated by a GP with zero mean and 
	covariance function $\kappa((\theta,\nu),(\tilde \theta,\tilde \nu))$, the posterior distribution of    $\JC(y_{1:\TT},u_{1:\TT};\theta^\star,\nu^\star)$ for generic controller parameters $(\theta^\star,\nu^\star)$ can be computed analytically. Specifically, $\JC(y_{1:\TT},u_{1:\TT};\theta^\star,\nu^\star)$ is a Gaussian   variable with mean 
	\begin{subequations} \label{eqn:GPmeanS}
		\begin{align}
		m_{i}(\theta^\star,\nu^\star) = & \bold{k}_{i}'\left(\bold{K}_{i} + \sigma_e^2I \right)^{-1}\JC_{1:i}, 
		\end{align} 
		and variance
		\begin{align}
		\sigma_{i}^2(\theta^\star,\nu^\star) = & \kappa((\theta^\star,\nu^\star),(\theta^\star,\nu^\star)) + \\
		- &   \bold{k}_{i}'\left(\bold{K}_{i} + \sigma_e^2I \right)^{-1}\bold{k}_{i} + \sigma_e^2,
		\end{align}	
	\end{subequations} 
	where the   $j$-th element of the vector $\bold{k}_{i}  \in\mathbb{R}^{i}$ is $\kappa((\theta^\star\!,\!\nu^\star),(\theta_j,\!\nu_j))$;    the $[j,h]$-th entry of the Kernel matrix  $\bold{K}_{i} \in \mathbb{R}^{i \times i}$ is $\kappa((\theta_j,\nu_j),(\theta_h,\nu_h))$;  $\sigma_e^2$ represents the variance of an additive (Gaussian) noise possibly affecting the observations of the cost $\JC$; and  $I$ denotes the identity matrix of proper dimension.
	
	The covariance function $\kappa((\theta,\nu),(\tilde \theta,\tilde \nu))$ for the GP can be defined, for instance, in terms of  
	the so-called \emph{Squared Exponential} (SE) covariance
	kernel, defined as 
	\begin{align*}
	\kappa((\theta,\nu),(\tilde \theta,\tilde \nu)) 
	= \sigma_0^2 e^{-\frac{1}{2\lambda^2} \left[\theta' - \tilde \theta' \ \mu' - \tilde \mu' \right] \left[\theta' - \tilde \theta' \ \mu' - \tilde \mu' \right]' }.
	\end{align*}
	The hyper-parameters  $\sigma_0$ and $\lambda$ characterizing the SE kernel, as well as the noise variance $\sigma_e^2$, can be chosen by  maximizing the
	log marginal likelihood \cite{williams2006gaussian}
	\begin{align*}
	& \log p(\JC_{1:i}|\theta_{1:i},\nu_{1:i},\sigma_0,\lambda,\sigma_e) \propto \\
	\propto & -\frac{1}{2}\log \det\left(\bold{K}_{i} + \sigma_e^2I \right)  - \frac{1}{2}\JC_{1:i}'\left(\bold{K}_{i} + \sigma_e^2I \right)^{-1}\JC_{1:i}. 
	\end{align*}
	
	%
	%
	\item \emph{Optimization phase} (Steps~\ref{algo:BOv2.2}-\ref{algo:BOv3}). In this phase, the next design parameters $\theta_{i+1}$ and $\nu_{i+1}$ to test are chosen by maximizing the so-called  \emph{acquisition function} $\alpha(\theta,\nu|\mathcal{D})$ (Step~\ref{algo:BOv2.2}). The acquisition function $\alpha(\theta,\nu|\mathcal{D})$ is constructed based on the mean and covariance (eq.~\eqref{eqn:GPmeanS}) of the GP estimated in the learning step. The acquisition function balances  exploration (\emph{i.e.}, learning more about the objective $\JC$ in regions of the parameter space with high variance) and exploitation (\emph{i.e.}, search over regions with high mean to optimize the expected performance based on past collected data). Different acquisition functions have been proposed in the literature (see \cite{brochu2010tutorial} and the references therein for a deep overview). In the example reported in Section~\ref{Sec:example}, we use the \emph{Expected Improvement} (EI) acquisition function, defined as 
	 \begin{align} \label{eqn:EI0}
	\alpha(\theta,\nu|\mathcal{D})=\textbf{EI}(\theta,\nu) = \mathbb{E}\left[\max\{0,\JC^- \!\!-\! \JC(\theta,\nu)\} \right],
	\end{align}
	where $\JC^-$ represents the best value of objective function at the $i$-th iteration, \emph{i.e.},
	\begin{align}
	\JC^- = \min\limits_{j=1,\ldots,i} \JC(y_{1:\TT},u_{1:\TT};\theta_j,\nu_j).
	\end{align}
	Under the GP framework previously discussed, the EI in~\eqref{eqn:EI0}  can be evaluated analytically and it is equal to:
	\begin{align}  \label{eqn:EI1}
	\textbf{EI}(\theta,\nu)\! =\!\!	\left(\!\JC^-\! - m_{i}(\theta,\nu)\right)\!\Psi(Z)\! +\!\sigma_{i}(\theta,\nu)\psi(Z)  
	\end{align}
	if $\sigma_{i}(\theta,\nu)>0$, $0$ otherwise. In \eqref{eqn:EI1}, $Z$ is defined as 
	\begin{align}
	Z = \frac{ \JC^- - m_{i}(\theta,\nu)  }{\sigma_{i}(\theta, \nu)}, 
	\end{align} 
	and  $\psi$ and $\Psi$ are the \emph{probability density function} and the \emph{cumulative density function} of the standard normal distribution, respectively. 
\end{itemize}

The advantages of using BO for tackling this design problem  are twofold. First, it is a derivative-free optimization algorithm, which is useful since a closed-form expression of the performance $\tilde J$ as a function of the design parameters $\theta$ and $\nu$ is not available. Second, it allows us to tune the controller parameters with as few evaluations of $\tilde J$ as possible. The latter point is crucial,  since each evaluation can be costly and time-consuming, as it requires a closed-loop experiment.

\subsection{Restricting the parameter space}

Bayesian optimization allows setting bounds on the  search space of the parameters $\theta$ and $\nu$. These bounds can be included in the maximization of the acquisition function at Step 2.3 of Algorithm 1.
Restricting the search space  generally speeds up the algorithm's convergence, thus requiring fewer evaluations of the functional $\tilde J$. 
Suitable bounds may be defined exploiting prior system knowledge and design choices. Some applicable 
restrictions of the parameter space are discussed next.

It may be reasonable to assume that the optimal solution is achieved using an inner controller $\mathcal{K}$ that stabilizes the inner loop $\mathcal{M}$.  Therefore, one may constrain $\mu$ so that the prediction model $M$ used by MPC is asymptotically stable. 

Some basic control design rules may be also used to restrict the search space of $\theta$ defining the inner-loop controller $\mathcal{K}$.
For instance,  if $\mathcal{K}$ is a PID controller parametrized as in~\eqref{eq:PID}, its static gain should have the same sign of the static gain of the (stable) system $\mathcal{S}$.

Another significant reduction of the parameter space may be achieved under the assumption that the prediction sub-model $M_y(\mu_y)$ used by the MPC accurately describes the system dynamics $\mathcal{M}$ from $g$ to $y$. In this case, one can simply derive the augmented model $M$ providing the relationship from $g$ to the plant input $u$ and output $y$  
as
\begin{equation}
\label{eq:aug_model_trick}
\matrice{c}{u \\ y} =  \underbrace{\matrice{c}{K(\theta)(I-M_y(\mu_y)) \\ M_y(\mu_y)}}_{M(\mu_y,\theta)}  g.
\end{equation}
Note that in this case $\mu=\smallmat{\mu_y\\\theta}$, that is the prediction model
and controller share some parameters.

Other restrictions may be introduced according to the particular problem at hand and prior knowledge  available to the user, \emph{e.g.}, diagonal models assuming decoupled
dynamics, grey-box models with known intervals for physical parameters, etc.

\section{Numerical Example}\label{sec:example}
\label{Sec:example}
\begin{figure}[!b]
	\centering
	\includegraphics[width=90pt]{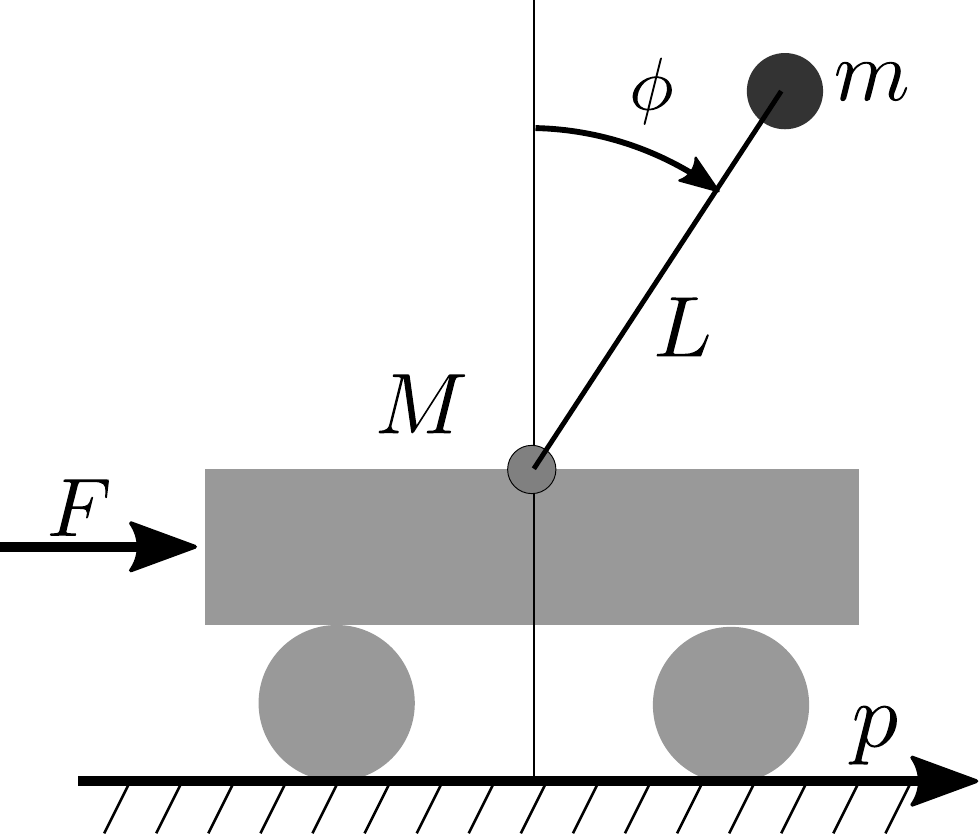} 
	\caption{Schematics of the inverted pendulum on a cart.}
	\label{fig:pendulum_cart}  
\end{figure}
As a case study, we consider the control problem of the inverted pendulum on a cart depicted in Fig.~\ref{fig:pendulum_cart}.  

\subsection{System description}
The dynamics of the process are governed by the  equations
\begin{subequations} \label{eqn:pend}
	\begin{align} 
	(M+m)\ddot p + mL\ddot\phi \cos\phi - mL \dot \phi ^2 \sin \phi + b\dot p &= F, \\
	L \ddot \phi + \ddot p \cos \phi - g \sin\phi + f_\phi\dot \phi &=0 ,
	\end{align}
\end{subequations} 
where $p$ is the cart position, $\phi$ is the angle of the pendulum with respect to the upright vertical position, and $F$ is an input force acting on the cart. The following values of the physical parameters are used:  $M=0.5~\text{Kg}$ (cart mass), $m=0.2~\text{Kg}$ (pendulum mass), $L=0.3~\text{m}$ (rod length), 
$g=9.81~\text{m/s\textsuperscript{2}}$ (gravitational acceleration), $b=0.1~\text{N/m/s}$, and $f_\phi=0.1~\text{m/s}$ (friction terms). According to the approach proposed in the paper, no knowledge of the physical model of the process is used in  designing the controller, and~\eqref{eqn:pend} are only used for data generation and performance evaluation. 



The output signals $p$ and $\phi$
are measured every $\Ts=5\;\text{ms}$ and measurements are  corrupted by an additive zero-mean white Gaussian noise with standard deviation $0.01~\text{m}$ and $0.01~\text{rad}$, respectively. The input force $F$ is also  perturbed by an additive zero-mean random disturbance with standard deviation $1~\text{N}$ and bandwidth $10\;\text{rad/s}$. 

In performing  closed-loop experiments, the system is initialized at $[p(0) \;  \dot{p}(0)\; \phi(0)\;\dot{\phi}(0)] = [0 \;  0\; \frac{\pi}{20}\;0]$. The   objective is to move the pendulum to the vertical position $\phi=0$, while limiting the cart displacement. 
The force $F$ is constrained to belong to the interval  $I_F=[-20\; 20]$~N, while the cart position $p$ should stay within the range $I_p=[-1\   1]$~\text{m} (representing, \emph{e.g.}, finite length of the track where the cart moves).


\subsection{Control design}
The hierarchical controller in Fig.~\ref{fig:hierarchical_outinner} is designed, with $y=[p\ \  \phi]$ and $u=F$.
The inner-loop controller $\Kin$ is
\begin{equation} \label{eqn:exKin}
u = \left[0\; \; \ \ \mathcal{K}_{PI}(z,\theta) \right](\rK-y),
\end{equation}
where $\mathcal{K}_{PI}(z,\theta)$ is a discrete-time transfer function of a PID controller parametrized as in \eqref{eq:PID},
with $\theta = [\theta_P\  \theta_I\  \theta_D] \in \mathbb R^3$ and $N_d=100$. Note that only the angle $\phi$ is actually  fed back in the inner loop, thus the task of the inner controller $\mathcal{K}$ is only to stabilize the dynamics of the angle $\phi$.

Besides taking care of  the   control objectives, the outer MPC shall also enforce    constraints on the cart position $p$ and on the input force $F$. 
The   model  used by the MPC to predict the dynamics of the inner loop $\mathcal{M}$ from the MPC command $\rK$ to the plant output  $y=[p \ \ \phi]$  (see Fig.~\ref{fig:hierarchical_outinner}) is parameterized as the continuous-time state-space model
\begin{equation} \label{eqn:exMM}
\dot \xi_M = A_{M}(\mu) \xi_M + B_{M}(\mu) \rK, \ \ \ y_M    = \xi_M, 
\end{equation}
where $\xi_M\in\mathbb{R}^2$ is the state vector and $\mu \in \mathbb R^{6}$ contains the  entries of   $A_M \in \mathbb R ^{2\times 2}$ and the second column of $B_M \in \mathbb R ^{2 \times 2}$. Because of the structure of the inner controller $\Kin$ in~\eqref{eqn:exKin}, the position $p$ is not fed back to the inner loop. Thus,  the first column of $B_M$ is set to zero and   not included in the design  parameter vector $\mu$. 
The overall MPC prediction model $M$ 
is constructed using~\eqref{eq:aug_model_trick}. 

The MPC control law is computed solving~\eqref{eq:MPC} and applied in a receding-horizon fashion, using a sampled version of~\eqref{eqn:exMM} 
with sampling time  $\TMPC = 10 \Ts = 50\;\text{ms}$,    
reference $r = [r_p \  r_\phi] = [0 \ \  0]$ and real-time constraints on $F$ and $p$ based on the  admissible intervals $I_F$ and $I_p$, respectively. 

Regarding the MPC design parameters, the weight matrices are not optimized and set  to 
$Q_y=\text{diag}(0.1,0.1)$, $Q_u = 0$,  $Q_{\Delta u} = 0.1$ and $Q_\epsilon = 10^5$. The prediction horizon $N_p$ is considered as a free parameter to be adjusted in the Bayesian optimization, while $N_u$ is set equal to $N_p$. The real-value design parameters $\theta$ and $\mu$ are constrained to belong to the interval $[-500 \ \ 500]$, while the prediction horizon $\Np$ can take integer values between $10$ and $20$. 
The MPC control law is computed using the MATLAB \emph{Model Predictive Control Toolbox}. All the computations are carried out on an i5 2.60-GHz Intel  core  processor  with  32  GB  of  RAM  running  MATLAB R2018a. 
The maximum computational time required to evaluate the MPC law over all the performed closed-loop experiments is  $21$ ms, thus lower than the  sampling time $\TMPC=50$ ms. 


Overall, there are $10$   parameters to be designed, namely $\theta \in \mathbb{R}^3$, $\mu \in  \mathbb{R}^6$, and $\Np \in  \mathbb{N}$. The  closed-loop performance cost $\JC$ to be minimized is defined as 
\begin{multline} \label{eqn:JCex}
\JC = \log \bigg[\frac{1}{T} \sum_{t=1}^T \bigg (\frac{1}{10}|r_p - p(t)|  + \frac{9}{10}|r_\phi - \phi(t)| \bigg ) \bigg ] +\\
+ \log \bigg[ \frac{1}{T}\sum_{t=1}^T b(p(t)) \;\;+1\bigg],
\end{multline}
where 
\begin{equation}
b(p) = \begin{cases}
10(|p|-1) \;\;\;  \;\;\;\text{for } |p| > 1,\\
0 \qquad \qquad \;\;\;\;\;\   \text{for } |p| \leq  1.
\end{cases}
\end{equation}
is the barrier function taking into account violation on the physical constraints on the cart position $p$.
The   cost $\JC$   is evaluated over closed-loop  experiments of length $10$~s on the discrete-time samples collected at rate $T_s$. This objective function reflects the engineering objective of controlling the angle $\phi$ to 0, limiting the horizontal displacement and keeping the cart position in the admissible range $I_p$. The constraint on the force $F$ is enforced by a saturation block at the system input, and thus it is not penalized  in $\JC$. 

The design problem \eqref{eqn:costCLtot} is solved using the
MATLAB \emph{Statistics and Machine Learning Toolbox}, setting the \textbf{EI} in \eqref{eqn:EI0} as acquisition function. $\Nin = 10$ random  values of the design parameters $\theta$, $\mu$ and $\Np$ are generated  to initialize Algorithm \ref{algo:BO}, which is then executed for $310$ iterations.  
The complete test code of this paper is available for download at \url{http://www.marcoforgione.it/data/code/CSL2019_perf.zip}.

\subsection{Simulation results}

\begin{figure}[!b]
	\centering
	\includegraphics[width=174pt]{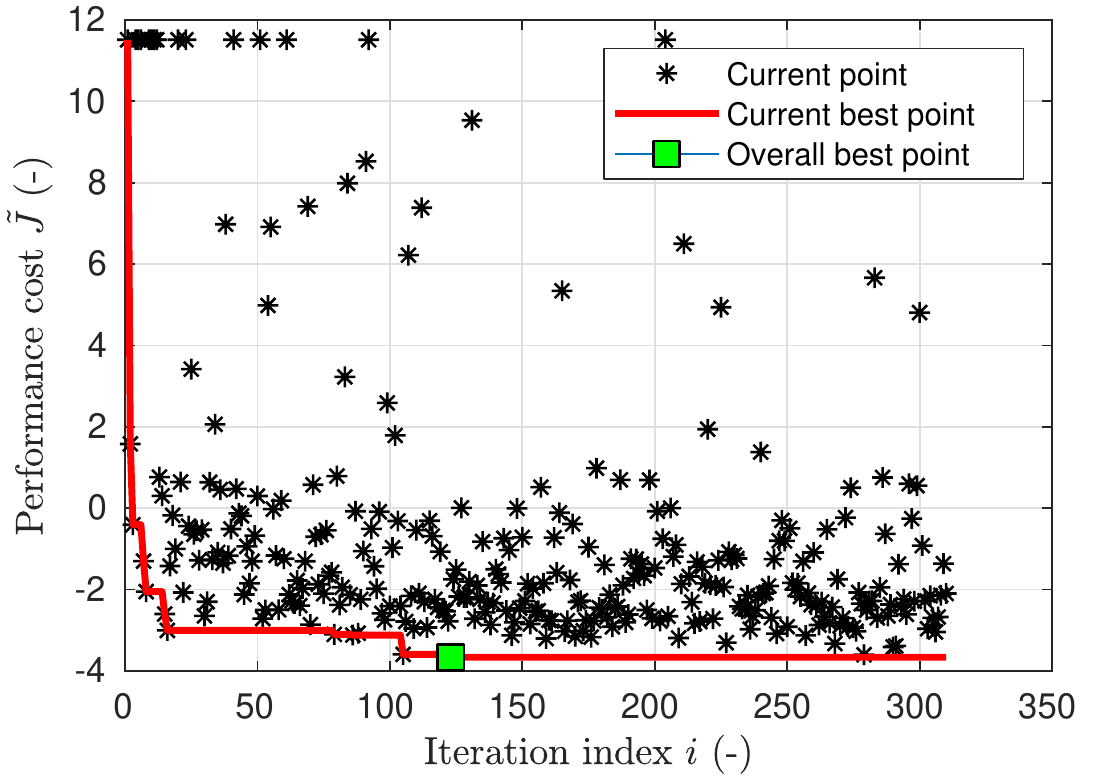}  
	\caption{Performance cost $\tilde J$ \emph{vs} iteration   $i$ of Algorithm~\ref{algo:BO}.}  
	\label{fig:iteration}
\end{figure}

The performance cost $\tilde J$ \textit{vs} the iteration index $i$ of Algorithm~\ref{algo:BO} is shown in Fig.~\ref{fig:iteration}. For each iteration $i$, the performance of the current test point (black asterisk) and of the current best point up to iteration $i$ (red line) are shown. 
From Fig.~\ref{fig:iteration}, it can be noticed that the  optimal controller parameters are found at iteration $123$ (green square). Furthermore,  as the iteration index $i$ increases, more and more test points are concentrated
in an area of low cost $\tilde J$. 

\begin{figure}[!b]
	\centering
	\includegraphics[width=174pt]{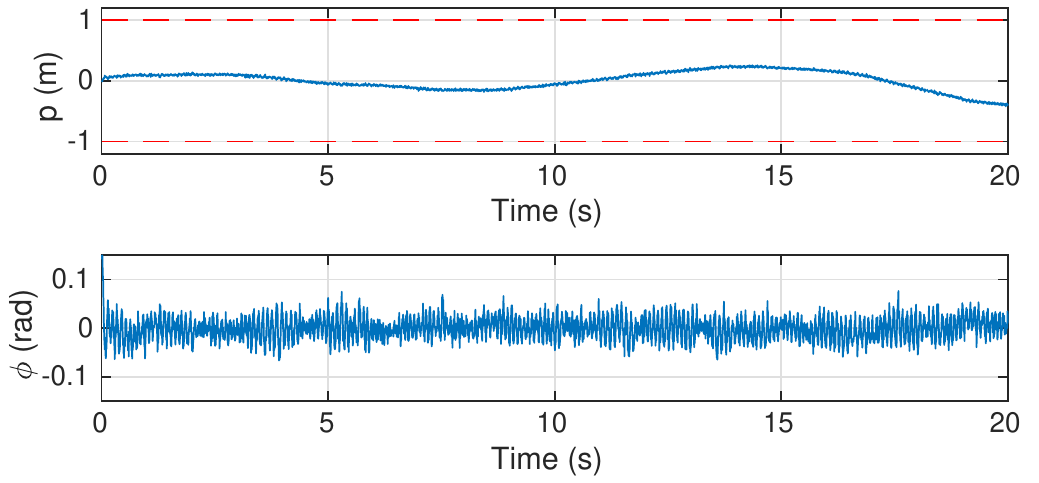}
 	\caption{Closed-loop experiment: position $p$ of the cart with the considered bounds (top plot) and pendulum's angle $\phi$ (bottom plot).}  
	\label{fig:position_angle}
\end{figure}

A closed-loop experiment is repeated over a longer period of $20$~s using the designed controller. The time trajectories of the cart's position $p$ and the pendulum's  
angle $\phi$ are plotted in Fig.~\ref{fig:position_angle}, which shows that the designed controller is able to stabilize the pendulum' angle in the upright vertical position, respecting the constraints on the cart's position $p$.  

For the sake of comparison, the following two non-hierarchical  model-based controllers are designed  based on the physical model of the system (Eq.~\eqref{eqn:pend}) linearized around   $[p(0) \;  \dot{p}(0)\; \phi(0)\;\dot{\phi}(0)] = [0 \;  0\; 0 \;0]$:
\begin{itemize}
	\item an MPC,  with the same sampling rate $\TMPC=50$~ms considered before, which reflects   real-time constraints.  At this sampling rate, the MPC is not able to reject the disturbance and thus fails to stabilize the pendulum around the upright vertical position.  This  shows the advantages of the  hierarchical multi-rate controller structure.
	
	\item a \emph{Linear-Quadratic-Gaussian} (LQG) controller, with sampling rate $\Ts=5$~ms. This controller stabilizes the pendulum. However, besides requiring a knowledge of the plant,  it achieves a performance cost $\JC$ (eq. \eqref{eqn:JCex}) equal to $-2.41$, which is worse than   the cost $\JC=-3.66$ obtained using  the proposed performance-oriented approach.  
\end{itemize} 


\section{Conclusions and follow-up}

In this work, we described a method to learn \textit{MPC-oriented models} for hierarchical control schemes via iterative closed-loop experiments. We showed that such experiments can be suitably designed using Bayesian optimization. In the proposed learning framework, the model does not necessarily provide the highest input/output data fit, which is the typical objective of system identification, but is the one yielding the model-based controller corresponding to the best closed-loop performance. We also argued that the prediction horizon can be optimized using the same tools and experiments. Numerical simulations on a benchmark example showed that data can lead to satisfactory controllers with no knowledge of the system dynamics and no constraints on modeling accuracy. Future research will be devoted to the theoretical analysis of the proposed learning strategy as well as to its experimental validation on a real-world setup.





\bibliographystyle{IEEEtran} 
\bibliography{Biblio_CDC}


\end{document}